\documentclass[12pt,letterpaper]{article}
\usepackage[utf8]{inputenc}
\usepackage[english]{babel}
\usepackage[T1]{fontenc}
\usepackage{amsfonts}
\usepackage{amsmath}
\usepackage{amstext}
\usepackage{amsbsy}
\usepackage{amssymb}
\usepackage{graphics}
\usepackage{graphicx}
\usepackage{yfonts}
\usepackage{pstricks-add}
\usepackage{pst-math}
\usepackage{url}

\usepackage{authblk}

\newtheorem{thm}{Theorem}[section]
\newcommand{\bthm}{\begin{thm}\noindent{\bf }~}
\newcommand{\ethm}{\end{thm}}
\newtheorem{propo}{Proposition}[section]
\newcommand{\bprop}{\begin{propo}\noindent{\bf }~}
\newcommand{\eprop}{\end{propo}}
\newtheorem{lema}{Lemma}[section]
\newcommand{\blem}{\begin{lema}\noindent{\bf }~}
\newcommand{\elem}{\end{lema}}
\newtheorem{defe}{Definition}[section]
\newcommand{\bdefe}{\begin{defe}\noindent{\bf }~}
\newcommand{\edefe}{\end{defe}}
\newtheorem{ejem}{Example}[section]
\newcommand{\bejem}{\begin{ejem}\noindent{\bf }~}
\newcommand{\eejem}{\end{ejem}}
\newtheorem{coro}{Corollary}[section]
\newcommand{\bcor}{\begin{coro}\noindent{\bf }~}
\newcommand{\ecor}{\end{coro}}
\newtheorem{rema}{Remark}[section]
\newcommand{\brem}{\begin{rema}\noindent{\bf }~\rm}
\newcommand{\erem}{\end{rema}}

\newcommand{\edemo}{\hfill{$\sqcap \kern-18pt \sqcup$}}

\title{\large\bf Potentials for solenoidal fields using the three-dimensional $\varphi-$harmonic cyclic algebra}

\author[1]{\small Homero G. D\'\i az-Mar\'\i n}
\author[2]{\small Elifalet L\'opez-Gonz\'alez}
\author[3]{\small Osvaldo Osuna}

\affil[1]{
		\footnotesize
		Facultad de Ciencias F\'\i sico-Matem\'aticas\\
		Universidad Michoacana\\
		Ciudad Universitaria, C.P. 58040 Morelia, M\'exico\\
		\ttfamily{homero.diaz@umich.mx}
}
\affil[2]{
      	 	\footnotesize
		DM de la UACJ en Cuauhtémoc\\
		Universidad Autónoma de Ciudad Juárez\\
       		Carretera Cuauhtémoc-Anáhuac, Km 3.5 S/N, Ejido 
		Cuauhtémoc, C.P. 31600, Cd. Cuauhtémoc, Chih, México\\
      		 \ttfamily{elgonzal@uacj.mx}
}
\affil[3]{
     		\footnotesize
		Instituto de F\'\i sica y Matem\'aticas, 							Universidad Michoacana,
		Ciudad Universitaria, C.P. 58040, Morelia, M\'exico.
		\ttfamily{osvaldo.osuna@umich.mx}
}  

\begin{document}

\maketitle

\begin{small}\noindent \textbf{Abstract.}
Given a PDE in \cite{LMT-2023} it is proposed a method for constructing solutions by considering an associative real algebra $\mathbb A$, and a suitable affine vector field
$\varphi$ with respect to which the components of all the
functions $\mathcal L\circ\varphi$ are solutions, where $\mathcal
L$ is differentiable in the sense of Lorch with respect to
$\mathbb A$. When we consider the 3D cyclic algebra and a suitable 3D affine map $\varphi$ we get families of solutions for the Laplace equation with three independent variables.
\end{small}

\noindent\textbf{Keyword}: \texttt{Functions of hypercomplex
variables, Laplace operator, Solenoidal vector fields, Calculus
over algebras, Differentiation theory.}

\noindent \textbf{MSC[2020]:} 
30G35,
35J05, 
35J47, 
35A25,  
58C20.


\section*{Introduction}

\noindent Harmonic vector fields $\mathbf V$ in $\mathbb R^3$ satisfy that their components are harmonic, i.e. 
\begin{equation}\label{eq:harmonic-vf}
	\Delta \, \mathbf V=0.
\end{equation}
\noindent where the components $V^i,i=1,2,3,$ satisfy Laplace's equation. In 2D and 3D Laplace's equations are respectively given
by
\begin{equation}\label{Laplace}
   u_{xx}+u_{yy}=0,\qquad
    u_{xx}+u_{yy}+u_{zz}=0.
\end{equation}
A special class of harmonic vector fields are lamellar or solenoidal vector fields, i.e. those that are incompressible and irrotational,
\begin{equation}\label{eq:solenoid-vf}
	\mathrm{div}\,\mathbf V=0,\qquad 
		\mathrm{curl}\,\mathbf V=0,
\end{equation}
We recall that the component of complex analytic functions are harmonic
functions in 2D, solving on a simply connected region, each harmonic function is
the real part of a complex analytic function.

Work has been carried out to build solutions for the
3D Laplace's equation and other PDEs of mathematical physics, by
using hypervariables; see \cite{Beckh}, \cite{Ket}, \cite{Ket1},
\cite{Miles}, \cite{Pog}, and \cite{Wag}. In these the
differentiability in the sense of Lorch has been used (or some
weaker differentiability using that of G\^ateaux), see also
\cite{Blum}, \cite{Lor}, \cite{War1}, \cite{War2}, and
\cite{Shef}. Expositions have recently been made on these topics
by J. S. Cook in \cite{JSC}, and by S. A. Plaksa in \cite{Pla}. Several of the given references have conditions of the type
that there exists a harmonic algebra with basis $\{e_1,e_2\}$ or $\{e_1,e_2,e_3\}$, for which $e_1^2+e_2^2=0$ in 2D or
\begin{equation}\label{eq:harmonic-id}
	e_1^2+e_2^2+e_3^2=0
\end{equation}
in 3D is satisfied in order to construct solutions
of the PDEs considered (\ref{Laplace}).

It is known that there is obstruction, described by Mel'nichenko, for the existence of real 3D algebras $\mathbb
A$ where the harmonic identity 	\eqref{eq:harmonic-id} holds true. Recently, complex algebras have been introduced to deal with this difficulty in \cite{Pla}.

On the other hand, for PDEs of the form
\begin{equation}\label{casi:laplace}
    Au_{xx}+Bu_{xy}+Cu_{yy}=0,
\end{equation}
a linear planar vector filed $\varphi$ and a 2D algebra $\mathbb
A$ are given, such that the components of the $\varphi\mathbb
A$-differentiable functions define a complete solution of
(\ref{casi:laplace}), see \cite{ELG3}. But this does not have a
similar result for the 3D case. The 3D version of the above result
gives the general form of harmonic functions
$$
u=\int_{-\pi}^{\pi}f(z+ix\cos s+iy\sin s, s)\,ds,
$$
where differentiations with respect to $x$, $y$, and $z$ under the
sign of integration can be done for the function $f$, see
\cite{Whit}.

The $\varphi\mathbb A$-differentiability of functions
$\mathcal F$ is introduced in \cite{LMT-2023}, where the
following definition is given: let $\varphi,\mathcal F$ be
 $n$-dimensional vector fields
which are differentiable in the usual sense on an open set $\mathcal U$, and
$\mathbb A$ a $n$-dimensional real algebra which is associative,
commutative, and has identity. If $F'_\varphi$ is a vector field such
that $d\mathcal F_p=F'_\varphi(p)d\varphi_p$ for all $p\in \mathcal U$, we
say that \emph{$\mathcal F$ is $\varphi\mathbb A$-differentiable}
and $F'_\varphi$ is its \emph{$\varphi\mathbb A$-derivative}. This
differentiability has associated its corresponding generalized
Cauchy-Riemann equations, see Section \ref{solutions}.

\noindent Recently, in \cite{ELG3} it is showed the components of
$\varphi\mathbb A$-differentiable functions define solutions for
PDEs; for each PDE of the form (\ref{casi:laplace}), and an affine
planar vector field $\varphi(x,y)=(ax+by,cx+dy)$ it is constructed
a two dimensional algebra $\mathbb A$ such that the components of
the second order $\varphi\mathbb A$-differentiable functions are
solutions of this PDE. By using the generalized Cauchy-Riemann
equations it is proved that every solution is a component of a
$\varphi\mathbb A$-differentiable function. So, a complete
solution is obtained. In particular, for the 2D Laplace's equation
given in left PDE at (\ref{Laplace}), if $Ac^{2}+Bcd+Cd^{2}\neq
0$, then $\mathbb A=\mathbb C$ if and only if
\[
Ac^{2}+Bcd+Cd^{2}=-(Aa^{2}+Bab+Cb^{2}),\qquad
2Aac+B(ad+bc)+2Cbd=0.
\]
In \cite{LMT-2023} the Cauchy problem defined by PDEs of the form
\eqref{casi:laplace} and conditions of the type
\begin{equation}\label{ccp}
    u(x,0)=\sum_{k=0}^\infty a_kx^k,\qquad u_y(x,0)=\sum_{k=0}^\infty
    b_kx^k,
\end{equation}
is solved. The solutions are expressed by power series with
respect to $\mathbb A$.

In this paper we consider PDEs with three independent
variables; the class of PDEs of the form
\begin{equation}\label{pde1}
    u_{xx}+u_{yy}+u_{zz}=0.
\end{equation}
\noindent If we consider a PDE as (\ref{pde1}), we look
for $\varphi$ like (\ref{varfi}), and an algebra $\mathbb A$ such
that the components of all the $\varphi\mathbb A$-differentiable
functions are solutions of the given PDE. The vector field
$\varphi$, and the algebra $\mathbb A$ are determined by a
solution of a system of three algebraic equations, as we described above.

For the method presented here, given a PDE and a vector
field, and then we look for an algebra, which is determined by a
solution of a system of three quadratic algebraic equations in six
variables. Also, we can give a PDE and an algebra, and then we
look for the vector field, which is determined by a solution of a
system of three quadratic algebraic equations in nine variables. Another
possible way is to consider a system of three quartic algebraic equations in
fifteen variables whose solutions determine the vector field and
the algebra.

The method applied in this article is a more explicit
way of that proposed in \cite{Ket} for solving PDEs of
mathematical physics, since here a more tractable type of
algebras, and the $\varphi\mathbb A$-differentiable functions are
used.

 In \cite{ELG3} it is used a family of
2D algebras
$$
\{\,\mathbb A^2_1(p_1,p_2)\,\,:\,\,p_1,\,p_2\in\mathbb R\,\}
$$
of two real parameters $p_1$, $p_2$, which are associative
commutative and have identity $e=e_1$, so that given a PDE from
mathematical physics (like the 2D Laplace's equation
\eqref{Laplace}), we look for a 2D affine transformation
$\varphi$, and an algebra $\mathbb A=\mathbb A^2_1(p_1,p_2)$ such
that condition of the type $\varphi(e_1)^2+\varphi(e_2)^2=0$ is
satisfied. In this work we use a family of six-parameter 3D
algebras
$$
\{\,\mathbb A^3_1(p_1,\cdots,p_6)\,\,:\,\,p_1,\cdots,p_6\in\mathbb
R\,\}
$$
which are associative, commutative, and have identity $e=e_1$, so that
given a PDE as \eqref{pde1}, we look for a transformation
$\varphi$ like (\ref{varfi}), and an algebra $\mathbb A=\mathbb
A^3_1(p_1,\cdots,p_6)$ such that the condition
\begin{equation}\label{vphihar3}
 \varphi(e_1)^2+\varphi(e_2)^2+\varphi(e_3)^2=0
\end{equation}
is satisfied. This is called a $\varphi$-harmonic algebra

\noindent In Section \ref{algebras} we introduce the algebras
considered in this paper. In Section \ref{solutions} we introduce
the $\varphi\mathbb A$-differentiability, and give a theorem about
solutions of PDEs with three independent variables. In Section
\ref{s3} we consider $\varphi\mathbb A$-harmonic algebras in 3D. We associate with each solution of a quartic system of
six algebraic equations in eighteen variables, an algebra and an
affine 3D vector field such that every $\varphi\mathbb
A$-differentiable functions has components which are solutions for
the 3D Laplace's equation. In Section \ref{s4} we obtain vector fields $\mathbf V$ solving \eqref{eq:solenoid-vf} from harmonic vector fields $\mathbf F$ solving \eqref{eq:harmonic-vf}.

\section{Pre-twisted real three dimensional algebras}

\subsection{Commutative algebras with identityy}\label{algebras}

\noindent We recall that a $\mathbb{R}$-linear space $\mathbb L$
is a \emph{commutative algebra with identityy} if it is endowed with a bilinear product
$\mathbb{L}\times\mathbb{L}\rightarrow\mathbb{L}$ denoted by
$(u,v)\mapsto u\centerdot v$, which is associative and commutative, $u\centerdot(v\centerdot w) =
(u\centerdot v)\centerdot w$ and $ u\centerdot v = v\centerdot u$ for all $u,v,w\in\mathbb{L}$; furthermore,
there exists an identity $e\in \mathbb{L},$ which satisfies $e\centerdot u=u$ for
all $u\in \mathbb{L}$. An element $u\in\mathbb L$ is called
\emph{regular} if there exists $u^{-1}\in\mathbb L$ called
\emph{the inverse} of $u$ such that $u^{-1}\centerdot u=e$. We also use the
notation $e/u$ for $u^{-1}$, where $e$ is the identity of $\mathbb L$.
If $u\in\mathbb L$ is not regular, then $u$ is called
\emph{singular}. $\mathbb L^*$ denotes the set of all the regular
elements of $\mathbb L$. If $u,v\in\mathbb L$ and $v$ is regular,
the quotient $u/v$ means $u\centerdot v^{-1}$.

\noindent It will be denoted by $\mathbb A$ if $\mathbb L=\mathbb
R^{3}$ and by $\mathbb M$ if $\mathbb L$ is a three dimensional
matrix algebra in the space of matrices $M(3,\mathbb R)$ where the
algebra product corresponds to the matrix product. We say that two
matrix algebras $\mathbb M_1$ and $\mathbb M_2$ are
\emph{conjugated} if there exists an invertible matrix $T$ such
that $\mathbb M_1=T\mathbb M_2T^{-1}$.

\noindent The $\mathbb A$ product between the elements of the
canonical basis $\{e_1,e_2,e_3\}$ of $\mathbb R^3$ is given by
\[
	e_i\centerdot e_j=\sum_{k=1}^3 c_{ijk}e_k
\]
where $c_{ijk}\in\mathbb R$ for
$i,j,k\in\{1,2,3\}$ are called \emph{structure constants} of
$\mathbb A$. The \emph{first fundamental representation} of
$\mathbb A$ is the injective linear homomorphism $R:\mathbb A\to
M(3,\mathbb R)$ defined by $R:e_i\mapsto R_i$, where $R_i$ is the
matrix with $[R_i]_{jk}=c_{ikj}$, for $i=1,2,3$.

Every three dimensional commutative algebra $\mathbb A$ with identity is isomorphic to one algebra belonging to a parametrized family {$\mathbb A^3_{r}({p_1},\cdots,{p_6})$ defined as follows.

\bdefe
The six parameter family of 3D algebras $\mathbb A^3_{r}({p_1},\cdots,{p_6})$ is the real linear space
$\mathbb R^{3}$ endowed with the product
\begin{equation}\label{algebra:3D:general}
 \begin{tabular}{c|ccc}
  $\centerdot$ & $e_r$ & $e_s$ & $e_t$ \\
  \hline
  $e_r$ & $e_r$ & $e_s$ & $e_t$ \\
  $e_s$ & $e_s$ & $p_7e_r+p_1e_s+p_2e_t$ & $p_8e_r+p_3e_s+p_4e_t$ \\
  $e_t$ & $e_t$ & $p_8e_r+p_3e_s+p_4e_t$ & $p_9e_r+p_5e_s+p_6e_t$ \\
\end{tabular},
\end{equation}
where the identities
\begin{equation}\label{ceq}
    \begin{array}{ccc}
      p_7 &=& -p_1p_4+p_2p_3-p_2p_6+p_4^2, \\
      p_8 &=& p_2p_5-p_3p_4,\qquad\qquad\quad\,\,\,\\
      p_9 &=& -p_1p_5+p_3^2-p_3p_6+p_4p_5, \\
    \end{array}
\end{equation}
stand for the associativity property. The identity is represented by $e=e_r$ in $\{e_r,e_s,e_t\}=\{e_1,e_2,e_3\}$. See
\cite{Dyg3} and \cite{Pie}. We recall that there are there are non-trivial isomorphisms classes as subsets of $\mathbb A^3_{r}({p_1},\cdots,{p_6})$.
\edefe

Moreover, the first fundamental representation
$R$ of $\mathbb A^{3}_{1}({p_1},\cdots,{p_6})$ is determined by
$$
R(e_1)=R_1=\left(%
\begin{array}{ccc}
  1 & 0 & 0 \\
  0 & 1 & 0 \\
  0 & 0 & 1 \\
\end{array}%
\right),\,\,\,
R(e_2)=R_2=\left(%
\begin{array}{ccc}
  0 & p_7 & p_8 \\
  1 & p_1 & p_3 \\
  0 & p_2 & p_4 \\
\end{array}%
\right),\,\,\,R(e_3)=R_3=\left(%
\begin{array}{ccc}
  0 & p_8 & p_9 \\
  0 & p_3 & p_5 \\
  1 & p_4 & p_6 \\
\end{array}%
\right).
$$
This allows us to use the corresponding matrix algebra in order to
get expressions of some vector fields which are defined with this
algebra product.

We will use extensively the {\em 3D cyclic algebra}, which corresponds to
\[
	\mathbb{A}= \mathbb A^{3}_1(0,1,0,0,1,0)
\]
which appears in
\cite{Ola} under the name of {\em Complex numbers in three dimensions} or {\em tricomplex numbers}. In \cite{Miles} it is used for constructing 3D harmonic functions.

The matrix algebra $\mathbb M=R(\mathbb A)$ is conjugated to the
matrix algebra spanned by the normal form with a real simple block
and a complex simple block, see \cite{AFLY} Section 2.2. Namely,
\begin{equation}\label{eq:rel}
R_1=
R_2^3=\left(%
\begin{array}{ccc}
  1 & 0 & 0 \\
  0 & 1 & 0 \\
  0 & 0 & 1 \\
\end{array}%
\right),\qquad
R_2=\left(%
\begin{array}{ccc}
  0 & 0 & 1 \\
  1 & 0 & 0 \\
  0 & 1 & 0 \\
\end{array}%
\right),\qquad R_3=R_2^2=\left(%
\begin{array}{ccc}
  0 & 1 & 0 \\
  0 & 0 & 1 \\
  1 & 0 & 0 \\
\end{array}%
\right).
\end{equation}

This cyclic algebra will be used in this paper for constructing
solutions for second order classical PDEs of the mathematical
physics: the 3D Laplace's equation.

Relations \eqref{eq:rel} for $e=e_1,e_2,e_3\in \mathbb A$ become
\begin{equation}\label{eq:rel2}
 \begin{tabular}{c|ccc}
  $\centerdot$ & $e$ & $e_2$ & $e_3$ \\
  \hline
  $e$ & $e$ & $e_2$ & $e_3$ \\
  $e_2$ & $e_2$ & $e_3$ & $e$ \\
  $e_3$ & $e_3$ & $e$ & $e_2$ \\
\end{tabular},
\end{equation}
where the identity matrix, $R(e)=R_1=\mathbb I_3$, corresponds to the identity, $e\in\mathbb A$.

We also remark that for the specific case of the cyclic algebra the set
\[
	\upsilon=x e+y e_2+z e_3\in\mathbb A
\]
is a regular element, i.e. $\upsilon\in \mathbb A^*$, if
\[
	\nu=x^3+y^3+z^3-3xyz=
	(x+y+z)(x^2+y^2+z^2-xy-yz-zx)\neq 0
\]
i.e. a so called {\em tricomplex number} $\upsilon\in\mathbb A$ has a unique inverse
\[
	\upsilon^{-1}=
	\frac{1}{\nu}\left[
	(x^2-yz)e
	+
	(z^2-xy)e_2
	+
	(y^2-zx)e_3\right],
\]
unless
\begin{equation}\label{eq:Pi}
x+y+z=0
\end{equation}
or
\begin{equation}\label{eq:trisector}
x^2+y^2+z^2-xy-yz-zx=0.
\end{equation}

We recall also the geometry of ${\mathbb A}^*$. According to \cite{Ola}, \eqref{eq:Pi} describes a plane, $\Pi\subset\mathbb A$, called {\em nodal plane}. On the other hand, relation \eqref{eq:trisector} becomes equivalent to the following condition for the so called {\em trisector line}, $\mathsf t\subset\mathbb A,$
\[
	x=y=z,
\]
which is perpendicular to $\Pi$ and generated by the vector
\begin{equation}\label{eq:n}
	\mathbf n=\frac{1}{\sqrt{3}}(e+e_1+e_2)\in\mathbb A,\qquad
	\mathbf n\cdot \Pi=0.
\end{equation}
Thus, the following assertion can be proved. 

\bprop \label{prop:1} The following properties hold true:
\begin{enumerate}
\item\label{inc:0} If $\upsilon'\in \Pi$, then $\upsilon\centerdot \upsilon' \in\Pi$ for every $\upsilon\in \mathbb A$.

\item\label{inc:1} If $\upsilon'\in \mathsf t$, then $\upsilon\centerdot \upsilon' \in\mathsf t$ for every $\upsilon\in \mathbb A$.

\item\label{inc:2} Whenever $\upsilon\in \Pi$ and $\upsilon'\in\mathsf t$, then
$
	\upsilon\centerdot \upsilon' = \mathbf 0.
$

\item\label{inc:2b} Whenever $\mathbf 0 \neq \upsilon'\in \Pi$ , $\upsilon\neq \mathbf 0$ and 
$
	\upsilon\centerdot \upsilon' = \mathbf 0,
$
 then $\upsilon\in\mathsf t$.

\item\label{inc:3} For any $\upsilon',\mu\in\Pi,$ $\upsilon'\neq\mathbf 0\neq \mu,$ there exists a solution $\omega\in \mathbb A$ of the equation 
\[
	\omega\centerdot \upsilon'=\mu.
\]
Whenever we look such solution conditioned to $\omega\in\Pi$, then such solution in unique.

\item\label{inc:4} If $\mathbf 0\neq \upsilon'\in \Pi$, and ${\upsilon}\in\mathbb A\setminus\mathsf t$, then there exists a unique $\omega\in \Pi$ such that
$
	\omega\centerdot \upsilon' =
	\upsilon\centerdot \upsilon' .
$
\end{enumerate}
\eprop

\noindent \emph{Proof.}
Straightforward calculations for $\upsilon=\upsilon^1e_1+\upsilon^2e_2+\upsilon ^3e_3$ and $\upsilon'=(\upsilon')^1e_1+(\upsilon')^2e_2+(\upsilon')^3e_3$ yield $\omega=\upsilon\centerdot\upsilon'$ as follows
\[\begin{aligned}
\omega^1&=
\upsilon^1(\upsilon')^1+\upsilon^2(\upsilon')^3+\upsilon^3(\upsilon')^2,
\\
\omega^2&=
\upsilon^1(\upsilon')^2+\upsilon^2(\upsilon')^1+\upsilon^3(\upsilon')^3,
\\
\omega^3&=
\upsilon^1(\upsilon')^3+\upsilon^2(\upsilon')^2+\upsilon^3(\upsilon')^1,
\end{aligned}\]
Thus, $(\upsilon')^1+(\upsilon')^2+(\upsilon')^3=0$ implies that $\omega^1+\omega^2+\omega^3=0$. This proves assertion \ref{inc:0}.

On the other hand, if $(\upsilon')^1=(\upsilon')^2=(\upsilon')^3$ then a simple inspection yields $\omega^1=\omega^2=\omega^3$. Thus we have proved claim \ref{inc:1}.

Since $\Pi\cap\mathsf t=\mathbf 0$, then properties \ref{inc:0} and \ref{inc:1} imply property \ref{inc:2}.

To prove \ref{inc:2b} we proceed as folles. First let us consider the following simplified expression of the product
\[
	\upsilon'\centerdot\upsilon
	=
	\left[
		(\upsilon')^1\upsilon^1+(\upsilon')^2\upsilon^3+(\upsilon')^3\upsilon^2
	\right]e+
	\left[
		(\upsilon')^1\upsilon^2+(\upsilon')^2\upsilon^1+(\upsilon')^3\upsilon^3
	\right]e_2+
	\left[
		(\upsilon')^1\upsilon^3+(\upsilon')^2\upsilon^2+(\upsilon')^3\upsilon^1
	\right]e_3.
\]
Then $\upsilon'\centerdot\upsilon=\mathbf 0$ becomes the homogeneous linear system
\[\begin{aligned}
	(\upsilon')^1\upsilon^1+(\upsilon')^2\upsilon^3+(\upsilon')^3\upsilon^2&=0,
	\\
	(\upsilon')^1\upsilon^2+(\upsilon')^2\upsilon^1+(\upsilon')^3\upsilon^3&=0,
	\\
	(\upsilon')^1\upsilon^3+(\upsilon')^2\upsilon^2+(\upsilon')^3\upsilon^1&=0,
\end{aligned}\]
which altogether with the orthogonality condition
\[
	(\upsilon')^1+(\upsilon')^2+(\upsilon')^3=0
\]
implies
\begin{equation}\label{system:2}
\begin{aligned}
	(\upsilon')^1(\upsilon^1-\upsilon^2)+(\upsilon')^2(\upsilon^3-\upsilon^2)&=0,
	\\
	(\upsilon')^1(\upsilon^2-\upsilon^3)+(\upsilon')^2(\upsilon^1-\upsilon^2)&=0.
\end{aligned}\end{equation}
If $(\upsilon')^1=0$ then \eqref{system:2} implies that $\upsilon^3=\upsilon^2=\upsilon^1$, whereas in the case $(\upsilon')^1=0$ we arrive at the same conclusion.

If $(\upsilon')^1\neq 0 \neq  (\upsilon')^2$ then \eqref{system:2} is a non-degenerate homogenous system in the variables $\upsilon^1-\upsilon^2, \upsilon^3-\upsilon^2$. Therefore, in any case the same conclusion arises. Namely,
\[\upsilon^3=\upsilon^2=\upsilon^1,\]
or equivalently $\upsilon\in \mathsf t$.

Uniqueness in claim \ref{inc:3} follows from
 \ref{inc:2} and  \ref{inc:2b} used for two possible solutions as follows:
 \[
 	\omega_1\centerdot \upsilon'=\mu
	=
	\omega_2\centerdot \upsilon'
	\quad \Rightarrow\quad
	(\omega_1-\omega_2)\centerdot\upsilon'=\mathbf 0
	\quad \Rightarrow\quad
	\omega_1-\omega_2\in\Pi \cap\mathsf t
	\quad \Rightarrow\quad
	\omega_1-\omega_2=\mathbf 0.
\]
Proving existence in claim \ref{inc:3} is equivalent to finding a unique solution $(\omega^1,\omega^2,\omega^3)$ of the linear system
\begin{equation}\label{eq:linear}\begin{aligned}
	\omega^1(\upsilon')^1+\omega^2(\upsilon')^3+\omega^3(\upsilon')^2
	&=\mu^1
	\\
	\omega^1(\upsilon')^2+\omega^2(\upsilon')^1+\omega^3(\upsilon')^3
	&=\mu^2
	\\
	\omega^1(\upsilon')^3+\omega^2(\upsilon')^2+\omega^3(\upsilon')^1
	&=\mu^3
	\\
	\omega^1+\omega^2+\omega^3&=0,
\end{aligned}\end{equation}
for fixed $(\mu^1,\mu^2,\mu^3)$ and $((\upsilon')^1,(\upsilon')^2,(\upsilon')^3)$ such that
\[
	\mu^1+\mu^2+\mu^3=0,
	\qquad
	(\upsilon')^1+(\upsilon')^2+(\upsilon')^3=0.                                   
\]
From linear dependence of coefficients, \eqref{eq:linear} can be reduced to
\[\begin{aligned}
	\omega^1(\upsilon')^1+\omega^2(\upsilon')^3+\omega^3(\upsilon')^2
	&=\mu^1
	\\
	\omega^1(\upsilon')^2+\omega^2(\upsilon')^1+\omega^3(\upsilon')^3
	&=\mu^2
	\\
	\omega^1+\omega^2+\omega^3&=0.
\end{aligned}\]
whose coefficients matrix has rank 3 for $\upsilon\neq \mathbf 0$ orthogonal to $(1,1,1)$.

Assertion \ref{inc:4} follows from claims \ref{inc:0} and \ref{inc:3}, with $\mu=\upsilon\centerdot\upsilon'$. \hfill$\square$

\bcor
The nodal plane $\Pi$ is an ideal of $\mathbb A$.
\ecor

\blem\label{lem:1}
Let us consider the linear map  $V: \mathbb A\rightarrow \mathbb A$
\begin{equation}\label{eq:V(F)}
V( \upsilon)= ( \upsilon^3-\upsilon^2)e_1+
			( \upsilon^3- \upsilon^1)e_2+
			( \upsilon^2- \upsilon^1)e_3
\end{equation}
Let $\upsilon_n:=\upsilon \cdot \mathbf n $ be normal component of $\upsilon$ and let 
\[
	\upsilon_\tau:=\upsilon - \upsilon_n\mathbf n\in \Pi,
\]
be its tangential component. Then
\[
	 V(\upsilon_\tau)=
	V(\upsilon)
\]
and 
\[
	\ker V=\langle\mathbf n\rangle=\mathsf t.
\]
\elem

\noindent \emph{Proof. }
Form linearity $V(\upsilon)=V(\upsilon_\tau)+\upsilon_n V(\mathbf n)$. On the other hand, $V(\mathbf n)=\boldsymbol{0}$. \hfill$\square$

\brem
From orthogonality, $\upsilon_\tau\cdot\mathbf n=0$ for $\upsilon_\tau=\upsilon_\tau^1e+\upsilon_\tau^2e_2+\upsilon_\tau^3e_3\in \Pi$, then
\[\begin{aligned}
	V(\upsilon)&=
	V(\upsilon_\tau)
\\
	&=
( \upsilon^3-\upsilon^2)e_1+
	( \upsilon^3- \upsilon^1)e_2+
	( \upsilon^2- \upsilon^1)e_3
\\
	&=
( -\upsilon^1-\upsilon^2-\upsilon^2)e_1+
	( -\upsilon^1-\upsilon^2- \upsilon^1)e_2+
	( \upsilon^2- \upsilon^1)e_3	
\\
	&=
	(-2\upsilon_\tau^2-\upsilon_\tau^1)e_1-
	(2\upsilon_\tau^1+\upsilon_\tau^2)e_2+
	(\upsilon_\tau^2-\upsilon_\tau^1)e_3.
\end{aligned}\]
\erem

For the following basis of $\mathbb A$,
\begin{equation}\label{eq:vi}
\begin{aligned}
	 v_1&=\frac{e+e_2+e_3}{3},\\
	 v_2&=\frac{2e-e_2-e_3}{3},\\
	 v_3&=\frac{e_2-e_3}{\sqrt{3}}.
	 \end{aligned}
\end{equation}
$v_2$ and $v_3$ are orthogonal, i.e. $v_2\cdot v_3=0$ and their Euclidean norm satisfy
\[\|v_2\|=\|v_3\|=\sqrt{2/3}.\]
The following relations can be checked, by straightforward calculations:
\begin{equation}\label{algebra:3Dv}
 \begin{tabular}{c|ccc}
  $\centerdot$ & $v_1$ & $v_2$ & $v_3$ \\
  \hline
  $v_1$ & $v_1$ 		& $\mathbf 0$ 	& $\mathbf 0$ \\
  $v_2$ & $\mathbf 0$ 	& $ v_2$ 		& $v_3$ \\
  $v_3$ & $\mathbf 0$ 	& $v_3$ 	& $-v_2$ \\
\end{tabular},
\qquad e=v_1+v_2.
\end{equation}
We also can conclude from the multiplication table \eqref{algebra:3Dv} the inclusion of the complex numbers as a subalgebra of the cyclic algebra $\mathbb A$. See Fig. \ref{fig:1} to see the 3D geometry of $\mathbb A$.

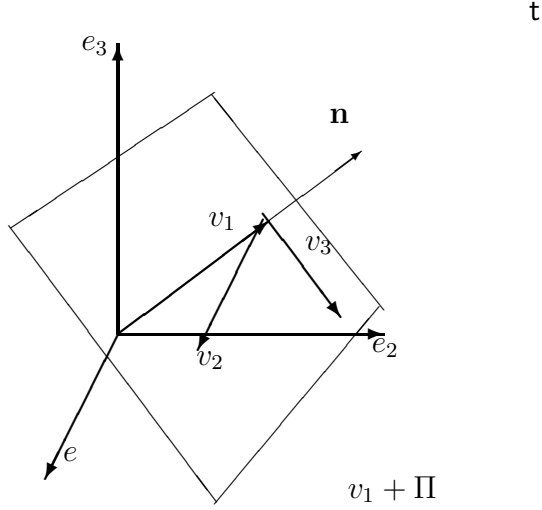
\begin{figure}[t] 
   \centering
   \setlength{\unitlength}{2000sp}%
\begingroup\makeatletter\ifx\SetFigFont\undefined%
\gdef\SetFigFont#1#2#3#4#5{%
  \reset@font\fontsize{#1}{#2pt}%
  \fontfamily{#3}\fontseries{#4}\fontshape{#5}%
  \selectfont}%
\fi\endgroup%
\begin{picture}(6849,6624)(3739,-8473)
\put(6076,-6736){\makebox(0,0)[lb]{\smash{{\SetFigFont{12}{14.4}{\rmdefault}{\mddefault}{\updefault}{\color[rgb]{0,0,0}$v_2$}%
}}}}
\thicklines
{\color[rgb]{0,0,0}\put(5101,-6361){\vector( 1, 0){3300}}
}%
{\color[rgb]{0,0,0}\put(5101,-6361){\vector( 0, 1){3600}}
}%
{\color[rgb]{0,0,0}\put(5101,-6361){\vector( 4, 3){1872}}
}%
{\color[rgb]{0,0,0}\put(6901,-4861){\vector( 3,-4){963}}
}%
\thinlines
{\color[rgb]{0,0,0}\put(5101,-6361){\vector( 4, 3){3024}}
}%
\thicklines
{\color[rgb]{0,0,0}\put(6901,-4936){\vector(-1,-2){810}}
}%
\thinlines
{\color[rgb]{0,0,0}\put(3751,-5011){\line( 3,-4){2574}}
\put(6301,-8461){\line( 5, 6){2040.984}}
\put(8401,-6061){\line(-4, 5){2136.585}}
\put(6301,-3361){\line(-3,-2){2526.923}}
}%
\put(4426,-7936){\makebox(0,0)[lb]{\smash{{\SetFigFont{12}{14.4}{\rmdefault}{\mddefault}{\updefault}{\color[rgb]{0,0,0}$e$}%
}}}}
\put(8251,-6586){\makebox(0,0)[lb]{\smash{{\SetFigFont{12}{14.4}{\rmdefault}{\mddefault}{\updefault}{\color[rgb]{0,0,0}$e_2$}%
}}}}
\put(4651,-2836){\makebox(0,0)[lb]{\smash{{\SetFigFont{12}{14.4}{\rmdefault}{\mddefault}{\updefault}{\color[rgb]{0,0,0}$e_3$}%
}}}}
\put(6226,-5011){\makebox(0,0)[lb]{\smash{{\SetFigFont{12}{14.4}{\rmdefault}{\mddefault}{\updefault}{\color[rgb]{0,0,0}$v_1$}%
}}}}
\put(7426,-5311){\makebox(0,0)[lb]{\smash{{\SetFigFont{12}{14.4}{\rmdefault}{\mddefault}{\updefault}{\color[rgb]{0,0,0}$v_3$}%
}}}}
\put(7726,-3736){\makebox(0,0)[lb]{\smash{{\SetFigFont{12}{14.4}{\rmdefault}{\mddefault}{\updefault}{\color[rgb]{0,0,0}$\mathbf{n}$}%
}}}}
\put(10201,-2461){\makebox(0,0)[lb]{\smash{{\SetFigFont{12}{14.4}{\rmdefault}{\mddefault}{\updefault}{\color[rgb]{0,0,0}$\mathsf{t}$}%
}}}}
\put(7951,-8386){\makebox(0,0)[lb]{\smash{{\SetFigFont{12}{14.4}{\rmdefault}{\mddefault}{\updefault}{\color[rgb]{0,0,0}$v_1+\Pi$}%
}}}}
\thicklines
{\color[rgb]{0,0,0}\put(5101,-6361){\vector(-1,-2){900}}
}%
\end{picture}   
   \caption{Geometry of $\mathbb A$}
   \label{fig:1}
\end{figure}

\bprop
The nodal plane constitutes a subalgebra
\[
	\Pi=\langle v_2,v_3\rangle =\langle \mathbf n\rangle^\bot \subset \mathbb A
\]
which is isomorphic to the complex numbers algebra, $\Pi \simeq \mathbb C,$ with isomorphism,
\[
	\Pi \ni av_2+bv_3 \leftrightarrow a+b\mathrm i\in\mathbb C,
	\qquad \forall a,b\in \mathbb R,
\]
given by the following identification:
$v_2\leftrightarrow 1,$ $ v_3\leftrightarrow {\rm i}, $ $v_3^2=-v_2\leftrightarrow  {\rm i}^2=-1.$ 
\eprop

\brem
Here we regard homomorphisms between the associative and commutative structure of the algebras regardless of the existence of identityy. Indeed, the identity $1\in\mathbb{C}$ corresponds to $v_2\in\Pi$ while $v_2\neq e\in\mathbb A$.
\erem

Since $V(\mathbf n)=\mathbf 0,$ the linear map $V:\mathbb A\rightarrow V(\mathbb A)$ has $\dim \ker  V=1$, then $V^{-1}(\beta_\star)$ has dimension 1 and is transverse to $\Pi$ for every $\beta_\star\in V(\mathbb A)$. Therefore, given any $\beta_\star\in V(\mathbb A)$, there exists a $\omega_\tau\in\Pi$ such that $V(\omega_\tau)=\beta_\star\in V(\mathbb A)$. Moreover, 
$
	V(\omega_\tau+\omega_n\mathbf n)=\beta_\star
$ for every $\omega_n\in\mathbb R$. Hence, without loss of generality we can suppose that $\omega_\tau\in\Pi$. Thus, we can prove the following Lemma.

\blem\label{lma:3}
The linear map $V\vert_\Pi: \Pi\rightarrow V(\Pi)$ is a linear isomorphism. Moreover,
\[\begin{aligned}
	V(\Pi)&= \left(e-e_2+e_3\right)^\bot,\\
	\ker V&=\langle\mathbf n\rangle,
\\
	V(\Pi)\cap\Pi
		&=
	\left\langle 
	V\left(v_2+\sqrt{3}v_3\right)\right\rangle=\langle e-e_3\rangle.
\end{aligned}\]
In particular $\Pi$ is not $V$-invariant.
\elem

\noindent \emph{Proof. }
We remark that the orthogonal basis $\{ v_2,v_3\}$ of $\Pi$, becomes the orthogonal basis $\{w_1,w_2\}=\{V(v_2),V(v_3)\}$ of $V(\Pi)$ given by
\[\begin{aligned}
	 w_2&=V(v_2)=2v_1-v_2=e_2+e_3,
	 \\
	 w_3&=V(v_3)=-\frac{2}{\sqrt{3}}v_1
	 		-\frac{2}{\sqrt{3}}v_2-v_3
			=\frac{1}{\sqrt{3}}\left[
				-2e-e_2+e_3
			\right],
\end{aligned}\]
which is also orthogonal with $
	\| w_2\|=\| w_3\|=\sqrt{2}.
$
\hfill$\square$

The geometry of $\Pi$ and $V(\Pi)$ described in Lemma \ref{lma:3} is illustrated in Fig.~\ref{fig:2}.

\begin{figure}[t] 
   \centering
   \setlength{\unitlength}{2000sp}%
\begingroup\makeatletter\ifx\SetFigFont\undefined%
\gdef\SetFigFont#1#2#3#4#5{%
  \reset@font\fontsize{#1}{#2pt}%
  \fontfamily{#3}\fontseries{#4}\fontshape{#5}%
  \selectfont}%
\fi\endgroup%
\begin{picture}(4749,5883)(3664,-8473)
\put(4051,-3586){\makebox(0,0)[lb]{\smash{{\SetFigFont{12}{14.4}{\rmdefault}{\mddefault}{\updefault}{\color[rgb]{0,0,0}$e-e_2+e_3$}%
}}}}
\thinlines
{\color[rgb]{0,0,0}\put(5626,-2911){\line( 2,-1){1980}}
\put(7576,-3961){\line(-1,-3){1492.500}}
\put(6151,-8461){\line(-1, 0){2475}}
\put(3676,-8311){\line( 1, 3){1815}}
}%
{\color[rgb]{0,0,0}\multiput(6301,-3361)(0.00000,-120.37037){41}{\line( 0,-1){ 60.185}}
}%
\thicklines
{\color[rgb]{0,0,0}\put(6226,-5686){\vector( 3, 2){1713.461}}
}%
{\color[rgb]{0,0,0}\put(6226,-5686){\vector(-1, 2){1035}}
}%
\put(7951,-8386){\makebox(0,0)[lb]{\smash{{\SetFigFont{12}{14.4}{\rmdefault}{\mddefault}{\updefault}{\color[rgb]{0,0,0}$\Pi$}%
}}}}
\put(8176,-4711){\makebox(0,0)[lb]{\smash{{\SetFigFont{12}{14.4}{\rmdefault}{\mddefault}{\updefault}{\color[rgb]{0,0,0}$\mathbf n$}%
}}}}
\put(5776,-2761){\makebox(0,0)[lb]{\smash{{\SetFigFont{12}{14.4}{\rmdefault}{\mddefault}{\updefault}{\color[rgb]{0,0,0}$V(\Pi)$}%
}}}}
\thinlines
{\color[rgb]{0,0,0}\put(3751,-5011){\line( 3,-4){2574}}
\put(6301,-8461){\line( 5, 6){2040.984}}
\put(8401,-6061){\line(-4, 5){2136.585}}
\put(6301,-3361){\line(-3,-2){2526.923}}
}%
\end{picture}%
   \caption{Planes $\Pi$ and $V(\Pi)$ inside ${\mathbb R}^3$.}
   \label{fig:2}
\end{figure}
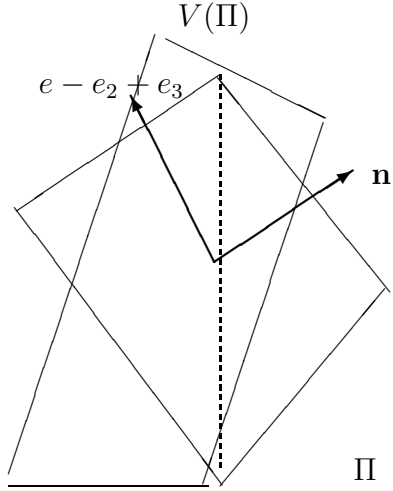

\subsection{Pre-twisted $\mathbb A$-differentiability and Cauchy-Riemann equations}\label{solutions}

The pre-twisted differentiability is introduced in \cite{LMT-2023},
this definition is closely related with the differentiability in
the sense of Lorch, see \cite{Lor}.

\begin{defe}
Let $\mathbb A$ be an algebra,
and
\[
	\varphi:\mathcal{U}\subset\mathbb{R}^3\rightarrow\mathbb A,
\]
a 3D vector field which is differentiable in the
usual sense. We say the 3D vector field, $\mathbf F:\mathcal{U}\subset\mathbb{R}^3\rightarrow\mathbb A,$ is
\emph{$\varphi\mathbb A$-differentiable (pre-twisted
differentiable)} if $\mathbf F$ is differentiable in the usual
sense and if there exists a 3D vector field $\mathbf F'_\varphi$
such that
\begin{equation}\label{dif:pre:algs}
    d\mathbf F_q=\mathbf F'_\varphi(q)\centerdot d\varphi_q,\qquad q \in \mathcal{U},
\end{equation}
 where $\mathbf F'_\varphi(q)\centerdot d\varphi_q(v)$ denotes the $\mathbb A$-product of $\mathbf F'_\varphi(q)$
 and $d\varphi(q)v$ for every vector $v$ in $\mathbb R^{3}$.
\end{defe}

A $\varphi\mathbb{A}$-\emph{polynomial function} $\mathbf
P:\mathbb R^3\to\mathbb A$ is defined by
\begin{equation}\label{Apolynomial:function}
    \mathbf
    P(q)=c_0+c_1\centerdot \varphi(q)+c_2\centerdot (\varphi(q))^{2}+\cdots+c_m\centerdot (\varphi(q))^m,
\end{equation}
where $c_0,c_1,\cdots,c_m\in\mathbb A$ are constants,
$q\in \mathcal{U}=\mathbb R^3$, and $c_k\centerdot(\varphi(q))^{k}$ for
$k\in\{1,2,\cdots,m\}$ are defined with respect to the $\mathbb
A$-product. If $\mathbf P$ and $\mathbf Q$ are
$\varphi\mathbb{A}$-polynomial functions, the
$\varphi\mathbb{A}$-\emph{rational function} $\mathbf P/\mathbf
Q$ is defined on the set $\mathbf Q^{-1}(\mathbb A^*)$. In the
same way \emph{exponential}, \emph{trigonometric}, and
\emph{hyperbolic $\varphi\mathbb A$-functions} are defined. All
these functions have $n$-order $\varphi\mathbb A$-derivatives for
$n\in\mathbb N$, and the usual rules for differentiation are
satisfied for this differentiability, except the chain rule.

The generalized Cauchy-Riemann equations for the differentiability
in the sense of Lorch can be seen in \cite{War2}. The
\emph{pre-twisted Cauchy-Riemann equations} associated with the
$\varphi\mathbb A$-differentiability were introduced in
\cite{LMT-2023}, they are given by the following relations among the first order partial derivatives,
 \begin{equation}\label{Pret-CR-EQs}
    \varphi_y \centerdot \mathbf F_x=\varphi_x \centerdot \mathbf F_y,\qquad \varphi_z \centerdot\mathbf
    F_x=\varphi_x \centerdot\mathbf F_z,\quad \varphi_z \centerdot\mathbf
    F_y=\varphi_y \centerdot\mathbf F_z.
\end{equation}
If $\varphi$ given in (\ref{varfi}) is an isomorphism, and
$\mathbf F$ is vector field which is differentiable in the usual
sense, then $\mathbf F$ is $\varphi\mathbb{A}$-differentiable if
and only if their components satisfy (\ref{Pret-CR-EQs}).

The first partial derivatives of every $\varphi\mathbb
A$-differentiable function $\mathbf F$ are expressed by
\begin{equation}\label{fpd:x:y:z}
\mathbf F_x = 
\mathbf F'_\varphi\centerdot\varphi_{x},\quad
\mathbf F_y =
\mathbf F'_\varphi\centerdot\varphi_{y},\quad
\mathbf F_z = \mathbf
F'_\varphi\centerdot\varphi_{z},
\end{equation}
while the second ones for an affine map $\varphi$ are given by
\begin{equation}\label{sp:xx:xy:yy:lineal}
\begin{array}{ccccccccccc}
             \mathbf F_{xx} & = & \mathbf F''_\varphi\centerdot\varphi_{x}^{2}, & \quad & \mathbf F_{yy} & = & \mathbf
             F''_\varphi\centerdot\varphi_{x}^{2}, & \quad & \mathbf F_{zz} & = & \mathbf F''_\varphi\centerdot\varphi_{z}^{2},\\
             \mathbf F_{xy} & = & \mathbf F''_\varphi\centerdot\varphi_{x}\centerdot\varphi_{y}, & \quad & \mathbf F_{xz} & = & \mathbf F''_\varphi\centerdot\varphi_{x}\centerdot\varphi_{z},
             & \quad & \mathbf F_{yz} & = & \mathbf F''_\varphi\centerdot\varphi_{y}\centerdot\varphi_{z}.
           \end{array}
\end{equation}

\subsection{Systems of algebraic equations associated with PDEs}

\noindent Given a PDE like (\ref{pde1}), we look for an affine change of coordinates as follows,
\begin{equation}\label{varfi}
   \varphi(x,y,z)=(a_1x+b_1y+c_1z+k_1,a_2x+b_2y+c_2z+k_2,a_3x+b_3y+c_3z+k_3)
\end{equation}
in a 3D algebra $\mathbb A$. From the product of $\mathbb A=\mathbb
A^{3}_1({p_1},\cdots,{p_6})$, and the proposed form for $\varphi$
in (\ref{varfi}) we have
\begin{equation}\label{dvarphi:x:cuadrado}
\begin{array}{ccc}
  \varphi_{x}^{2} & = & (a_1^{2}+a_2^{2}(-p_1p_4+p_2p_3-p_2p_6+p_4^{2})+a_3^{2}(-p_1p_5+p_3^{2}-p_3p_6+p_4p_5))e_1 \\
   &  & +2a_2a_3(p_2p_5-p_3p_4)e_1+(2a_1a_2+2a_2a_3+a_2^{2}p_1+a_3^{2}p_5)e_2 \\
   &   & +(2a_1a_3+a_2^{2}p_2+2a_2a_3p_4+a_3^{2}p_6)e_3,  \\
\end{array}
\end{equation}
\begin{equation}\label{dvarphi:y:cuadrado}
\begin{array}{ccc}
   \varphi_{y}^{2} & = & (b_1^{2}+b_2^{2}(-p_1p_4+p_2p_3-p_2p_6+p_4^{2})+b_3^{2}(-p_1p_5+p_3^{2}-p_3p_6+p_4p_5))e_1 \\
   &  & +2b_2b_3(p_2p_5-p_3p_4)e_1+(2b_1b_2+2b_2b_3+b_2^{2}p_1+b_3^{2}p_5)e_2 \\
   &   & +(2b_1b_3+b_2^{2}p_2+2b_2b_3p_4+b_3^{2}p_6)e_3,\\
\end{array}
\end{equation}
\begin{equation}\label{dvarphi:z:cuadrado}
\begin{array}{ccc}
   \varphi_{z}^{2} & = & (c_1^{2}+c_2^{2}(-p_1p_4+p_2p_3-p_2p_6+p_4^{2})+c_3^{2}(-p_1p_5+p_3^{2}-p_3p_6+p_4p_5))e_1 \\
   &  & +2c_2c_3(p_2p_5-p_3p_4)e_1+(2c_1c_2+2c_2c_3+c_2^{2}p_1+c_3^{2}p_5)e_2 \\
   &   & +(2c_1c_3+c_2^{2}p_2+2c_2c_3p_4+c_3^{2}p_6)e_3,  \\
\end{array}
\end{equation}
\begin{equation}\label{dvarphi:xy:cuadrado}
\begin{array}{ccc}
   \varphi_{x}\varphi_{y} & = &
   ({a_1}{b_1}+{a_2}{b_2}(-{p_1}{p_4}+{p_2}{p_3}-{p_2}{p_6}+{p_4}^2)
+{a_2}{b_3}(-{p_1}{p_4}+{p_2}{p_3}))e_1\\
   &  & +({a_2}{b_3}(-{p_2}{p_6}+{p_4}^2)+{a_3}{b_2}({p_2}{p_5}-{p_3}{p_4})+{a_3}{b_3}({p_2}{p_5}+{p_3}{p_4}))e_1 \\
      &  &  +({a_1}{b_2}+{a_1}{b_3}+{a_2}{b_1}+{a_2}{b_2}{p_1}+{a_2}{b_3}{p_1}+{a_3}{b_2}{p_3}+{a_3}{b_3}{p_3})e_2  \\
   &  &  +(a_3b_1+a_2b_2p_2+a_2b_3p_2+a_3b_2p_4+a_3b_3p_4)e_3,  \\
\end{array}
\end{equation}
\begin{equation}\label{dvarphi:xz:cuadrado} \begin{array}{ccc}
   \varphi_{x}\varphi_{z} & = &
   ({a_1}{c_1}+{a_2}{c_2}(-{p_1}{p_4}+{p_2}{p_3}-{p_2}{p_6}+{p_4}^2)
+{a_2}{c_3}(-{p_1}{p_4}+{p_2}{p_3}))e_1\\
   &  & +({a_2}{c_3}(-{p_2}{p_6}+{p_4}^2)+{a_3}{c_2}({p_2}{p_5}-{p_3}{p_4})+{a_3}{c_3}({p_2}{p_5}+{p_3}{p_4}))e_1 \\
      &  &  +({a_1}{c_2}+{a_1}{c_3}+{a_2}{c_1}+{a_2}{c_2}{p_1}+{a_2}{c_3}{p_1}+{a_3}{c_2}{p_3}+{a_3}{c_3}{p_3})e_2  \\
   &  &  +(a_3c_1+a_2c_2p_2+a_2c_3p_2+a_3c_2p_4+a_3c_3p_4)e_3,  \\
\end{array}
\end{equation}
\begin{equation}\label{dvarphi:yz:cuadrado} \begin{array}{ccc}
   \varphi_{y}\varphi_{z} & = &
({b_1}{c_1}+{b_2}{c_2}(-{p_1}{p_4}+{p_2}{p_3}-{p_2}{p_6}+{p_4}^2)
+{b_2}{c_3}(-{p_1}{p_4}+{p_2}{p_3}))e_1\\
   &  & +({b_2}{c_3}(-{p_2}{p_6}+{p_4}^2)+{b_3}{c_2}({p_2}{p_5}-{p_3}{p_4})+{b_3}{c_3}({p_2}{p_5}+{p_3}{p_4}))e_1 \\
      &  &  +({b_1}{c_2}+{b_1}{c_3}+{b_2}{c_1}+{b_2}{c_2}{p_1}+{b_2}{c_3}{p_1}+{b_3}{c_2}{p_3}+{b_3}{c_3}{p_3})e_2  \\
   &  &  +(b_3c_1+b_2c_2p_2+b_2c_3p_2+b_3c_2p_4+b_3c_3p_4)e_3.  \\
\end{array}
\end{equation}
If we consider that $a_i$, $b_i$, $c_i$ for $i=1,2,3$ as
fixed numbers, while $p_j$ for $j=1,\dots,6$ are variables,
then each of these equations corresponds to three quadratic
equations in six variables. On the other hand, if we consider that $a_i$, $b_i$,
$c_i$ for $i=1,2,3$ as variables, and  $p_j$ for
$j=1,\dots,6$ as fixed numbers, then each of these equations
corresponds to three quadratic equations in nine variables. Finally, if
we consider $a_i$, $b_i$, $c_i$ for $i=1,2,3$, as well as $p_j$ for
$j=1,\dots,6$ as variables, then each of these equations
corresponds to three quartic equations in fifteen variables.

\section{$\varphi$-harmonic algebras}\label{s3}

P. W. Ketchum called an algebra $\mathbb A$ a \emph{harmonic
algebra} if their analytic functions satisfy Laplace equation.
We introduce the following definition.

\bdefe If $\varphi$ is an affine vector field and $\mathbb A$ an
algebra such that the identity
\begin{equation}\label{par:varfi}
\begin{array}{c}
  \varphi_{x}^{2}+\varphi_{y}^{2}+\varphi_{z}^{2}=0, \\
\end{array}
\end{equation}
is satisfied, then
$\mathbb A$ will be called \emph{$\varphi$-harmonic algebra}.
\edefe

Note that $\varphi_x=d\varphi(e_1)$, $\varphi_y=d\varphi(e_2)$, and
$\varphi_z=d\varphi(e_3)$.

H. A. V. Beckh-Widmanstetter \cite{Beckh} has proved that there
does not exist a 3D harmonic algebra with identity $e=e_1$ over the
field $\mathbb R$. That is, there does not exist a 3D algebra $\mathbb
A$ with identity $e=e_1$ such that $e_1^2+e_2^2+e_3^2=0$.

On the contrary we provide conditions
for the $\varphi$-harmonicity of $\mathbb A=\mathbb
A^{3}_1(p_1,\cdots,p_6)$ which are given in the following proposition.

\bprop\label{con:alg:sist} Let $(p_1,\cdots,p_9)$ be a solution of
the system
\begin{equation}\label{cond:varphi:armonicas}
    \begin{array}{ccc}
  -x_1x_4+x_2x_3-x_2x_6+x_4^2-x_7 & = & 0,\\
  x_2x_5-x_3x_4-x_8 & = & 0,\\
  -x_1x_5+x_3^2-x_3x_6+x_4x_5-x_9 & = & 0,\\
  \|A_2\|^2x_1+2(A_2\cdot A_3)x_3+\|A_3\|^2x_5 & = & -2(A_1\cdot A_2), \\
  \|A_2\|^2x_2+2(A_2\cdot A_3)x_4+\|A_3\|^2x_6 & = & -2(A_1\cdot A_3), \\
  \|A_2\|^2x_7+2(A_2\cdot A_3)x_8+\|A_3\|^2x_9 & = & -\|A_1\|^2, \\
\end{array}
\end{equation}
where $A_i=(a_i,b_i,c_i)$, and $\cdot$ denotes the inner product
in $\mathbb R^3$. Thus, for $\mathbb A=\mathbb
A^3_1(p_1,\cdots,p_6)$ and $\varphi$ given by (\ref{varfi}),
$\mathbb A$ is $\varphi$-harmonic. \eprop \noindent\textbf{Proof.}
Let $(p_1,\cdots,p_9)$ be a solution of the system
(\ref{cond:varphi:armonicas}), and $\mathbb A=\mathbb
A^3_1(p_1,\cdots,p_6)$. By using the $\mathbb A$-product we obtain
\begin{equation}\label{Cond:harm}
    \begin{array}{ccc}
      \varphi(e_1)^2+\varphi(e_2)^2+\varphi(e_3)^2 & = & (\|A_1\|^2+\|A_2\|^2p_7+\|A_3\|^2p_9+2(A_2\cdot A_3)p_8)e_1 \\
       & + & (\|A_2\|^2p_1+\|A_3\|^2p_5+2(A_1\cdot A_2)+2(A_2\cdot A_3)p_3)e_2 \\
       & + & (\|A_2\|^2p_2+\|A_3\|^2p_6+2(A_1\cdot A_3)+2(A_2\cdot
    A_3)p_4)e_3.
    \end{array}
\end{equation}
From last three equations of system (\ref{cond:varphi:armonicas})
we obtain (\ref{vphihar3}). $\Box$

Proposition \ref{con:alg:sist} cannot be satisfied for orthonormal basis $\{A_1,A_2,A_3\}$ i.e. fro orthogonal matrix $A$.

It is satisfied for the 3D cyclic algebra as it is shown in the following assertion.

\bcor\label{harmonic1} For the affine map
\begin{equation}\label{varfi:p:acyc1}
    \varphi(x,y,z)=(-x-y+k_1,x-z+k_2,y+z+k_3),
\end{equation}
the algebra $\mathbb A^3_1(0,1,0,0,1,0)$ is a $\varphi$-harmonic
algebra. \ecor

\noindent\textbf{Proof.} For the algebra $\mathbb
A=\mathbb A^3_1(0,1,0,0,1,0)$ the parameters $p_i$ are given by
$p_1=0$, $p_2=1$, $p_3=0$, $p_4=0$, $p_5=1$, $p_6=0$, $p_7=0$,
$p_8=1$, and $p_9=0$.

For the vector field $\varphi$ given in (\ref{varfi:p:acyc1}) we
have that
$$
A_1 = (-1,-1,0),\quad A_2 =(1,0,-1),\quad A_3=(0,1,1).
$$
i.e. \begin{equation}\label{eq:phiA}
\varphi(q)=Aq+ k,\end{equation}
where
\begin{equation}\label{eq:A}
	A= 
	\left(\begin{array}{ccc}
		-1 & -1 & 0 \\
					   1 & 0 & -1 \\
					   0 & 1 & 1
	\end{array}\right)
	=
	\left(\begin{array}{c}
		A_1\\
		A_2\\
		A_3
	\end{array}\right)
	,
\end{equation}
$q=(x,y,z)^\dagger $ and $k=k^1e_1+k^2e_2+k^3e_3$. So that, $\|A_i\|^2=2$ for $i=1,2,3$, and
$$
A_1\cdot A_2=-1,\quad A_1\cdot A_3=-1, \quad A_2\cdot A_3=-1.
$$
Then, $p_i$ for $i=1,\cdots,9$ is a solution of system
(\ref{cond:varphi:armonicas}). Thus, by Proposition
\ref{con:alg:sist} we obtain (\ref{vphihar3}). Therefore, $\mathbb
A$ is a $\varphi$-harmonic algebra. $\Box$

\blem\label{lma:phi(A)} The linear map $A:\mathbb R^3\rightarrow\mathbb A$, induces an isomorphism $A\vert_\Pi:\Pi\rightarrow\Pi$. More precisely,
\[\begin{aligned}
\ker A&= \langle e-e_2+e_3
 \rangle,
 \\
 A(\Pi)&=\Pi=\langle\mathbf n\rangle^\bot.
 \end{aligned}\]
Thus, $\Pi\subset\mathbb R^3$ is $A$-invariant.
\elem

\noindent \emph{Proof. }
For the basis $\{e_2-e_3,e_2-e_1\}$ of the subspace $\Pi$, we have
\[\begin{aligned}
	A(e_2-e_3)=e_2-e_1,
	\qquad A(e_2-e_1)=-(e_2-e_3),
\end{aligned}\]
and $A(q)=\mathbf 0$ implies that $q^1=-q^2=q^3$. \hfill$\square$

Notice that $A\vert_\Pi$ is a $90^\circ$ rotation in $\Pi$, while
\[
	\varphi(e_2-e_3)=e_2-e_1+k,
	\qquad\varphi(e_2-e_1)=-(e_2-e_3)+k.
\]
The
 equalities (\ref{par:varfi}), give rise to
solutions of a PDEs. We can prove for instance the following assertion.

\blem
Let $\mathbb A$ be the cyclic 3D algebra, and $\varphi(q)=Aq+k$ be the
affine map (\ref{varfi}) where $A$ is the matrix \eqref{eq:A}. Then for a differentiable vector field, $\mathbf F:\mathcal U\rightarrow\mathbb A$, to be $\varphi\mathbb A$-differentiable is necessary and sufficient to satisfy four linearly independent Cauchy-Riemann (CR) such as:
\begin{equation}\label{eq:system2}
\begin{array}{rcl}
	F^1_x - F^1_y -F^2_x +  F_y^3& = &0,\\
	- F^1_x +  F_y^2+ F^3_x -  F^3_y & = & 0,\\
	- F^1_z - F^2_x +  F_x^3+ F^3_z & = & 0,\\
	 F^1_x+F^1_z    - F^2_z -  F_x^3& = & 0.
\end{array}
\end{equation}
\elem

\noindent\textbf{Proof.}
 The pre-twisted Cauchy-Riemann equations \eqref{Pret-CR-EQs} read as follows,
\begin{equation}\label{eq:CR-cyclic}
	b\centerdot \mathbf F_x - a\centerdot \mathbf F_y =0,
	\quad
	c\centerdot \mathbf F_x - a\centerdot \mathbf F_z =0,
	\quad
	c\centerdot \mathbf F_y - b\centerdot \mathbf F_z =0,
\end{equation}
where $a,b,c\in\mathbb A$ are the column vectors, $A=(a\,\vert\,b\,\vert c).$ For the cyclic algebra,
\[
	a=(-1,1,0)^\dagger,\quad
	b=(-1,0,1)^\dagger,\qquad
	c=(0,-1,1)^\dagger.
\]
More explicitly, we get a system of 9 homogeneous equations contained in \eqref{eq:CR-cyclic} as follows,
\[\begin{aligned}
	a_1 F^1_y - 	b_1 F^1_x +
	a_3 F^2_y - 	b_3 F^2_x +
	a_2 F^3_y - 	b_2 F^3_x
	&=0,\\
	a_2 F^1_y - 	b_2 F^1_x +
	a_1 F^2_y - 	b_1 F^2_x +
	a_3 F^3_y - 	b_3 F^3_x
	&=0,\\
	a_3 F^1_y - 	b_3 F^1_x +
	a_1 F^3_y - 	b_1 F^3_x +
	a_2 F^2_y - 	b_2 F^2_x
	&=0,\\
	\\
	a_1 F^1_z - 	c_1 F^1_x +
	a_3 F^2_z - 	c_3 F^2_x +
	a_2 F^3_z - 	c_2 F^3_x
	&=0,\\
	a_2 F^1_z - 	c_2 F^1_x +
	a_1 F^2_z - 	c_1 F^2_x +
	a_3 F^3_z - 	c_3 F^3_x
	&=0,\\
	a_3  F^1_z - 	c_3  F^1_x +
	a_1  F^3_z - 	c_1  F^3_x +
	a_2 F^2_z - 	c_2  F^2_x
	&=0,\\
	\\
	b_1  F^1_z - 	c_1  F^1_y +
	b_2  F^2_z - 	c_2  F^2_y +
	b_3  F^3_z - 	c_3  F^3_y
	&=0,\\
	b_2  F^1_z - 	c_2  F^1_y +
	b_1  F^2_z - 	c_1  F^2_y +
	b_3  F^3_z - 	c_3  F^3_y
	&=0,\\
	b_3  F^1_z - 	c_3  F^1_y +
	b_1  F^3_z - 	c_1  F^3_y +
	b_2 F^2_z - 	c_2  F^2_y
	&=0.	
\end{aligned}\]
For the specific case of the cyclic algebra, 
\begin{equation}\label{eq:system1}
\begin{array}{cc}
\begin{array}{rcl}
	- F^1_y +  F^1_x -F^2_x +  F_y^3& = &0,\\
	 F^1_y -  F^2_y + F^2_x - F_x^3& =& 0,\\
	- F^1_x -  F^3_y + F^3_x +  F_y^2& = & 0,
\end{array}
&
\begin{array}{rcl}
	- F^1_z - F^2_x + F^3_z +  F_x^3& = & 0,\\
	 F^1_z +  F^1_x - F^2_z -  F_x^3& = & 0,\\
	- F^1_x - F^3_z +F^2_z + F_x^2& = & 0,
\end{array}
\\
\\
\begin{array}{rcl}
	- F^1_z +  F^2_z -F^2_y +  F_y^3& = & 0,\\
	 F^1_y -  F^2_z + F^3_z - F_y^3& = & 0,\\
	F^1_z -  F^1_y -  F^3_z +  F_y^2& = & 0.
\end{array}
&
\end{array}
\end{equation}
Regarding \eqref{eq:system2} as a linear system of nine equations on nine variables, $F^{i}_a$ where $i=1,2,3,$ and $a\in\{x,y,z\}$, a straightforward calculation yields that it has rank 4. In addition, linear system \eqref{eq:system1} has also rank 4. Therefore, system \eqref{eq:system1} can be reduced to a system of four CR equations \eqref{eq:system2}. \hfill $\square$

\bthm\label{teo:0} Let $\mathbb A$ be the 3D cyclic algebra, and $\varphi(q)=Aq+k$ be the
affine map \eqref{eq:phiA} where $A$ is the matrix \eqref{eq:A}. Then the
components of a $\varphi\mathbb A$-differentiable vector field, $\mathbf F:\mathcal U\rightarrow\mathbb A$, satisfying the CR equations \eqref{eq:system2},
\[
	\mathbf F (x,y,z)= 
		F^1(x,y,z)e+
		F^2(x,y,z)e_2+
		F^3(x,y,z)e_3,
\]
are harmonic, i.e. the components $F^i$ are solutions of
(\ref{pde1}), and $\mathbf F$ solves the equation
\[
	\Delta\,\mathbf F=\mathbf 0.
\]
\ethm

\noindent\textbf{Proof.} 
If we multiply (\ref{par:varfi}) by
$\mathbf F''_\varphi$, and use (\ref{sp:xx:xy:yy:lineal}), we obtain that
components of $\mathbf F$ are solutions for (\ref{pde1}).
Explicitly, for the cyclic algebra and $k=\mathbf 0$ we have
\[
\varphi(e)^2+\varphi(e_2)^2+\varphi(e_3)^2=
(-e+e_2)^2+(-e+e_3)^2+(-e_2+e_3)^2=\mathbf 0.
\]
Thus, $\mathbf F$ is a $\varphi\mathbb A$-differentiable function that its components are solutions for (\ref{Laplace}). \hfill $\square$

\bejem
As we have mentioned we can consider any polynomial function
$
	\mathbf F(q)=c_0+\dots+c_n\varphi(q)^n,$ $q\in\mathbb R^3
$
which for $c_0\in\Pi$ is parallel to the plane $\Pi$, since $\Pi$ is an ideal of $\mathbb A$. Take for instance $\varphi(q)\in\Pi$ with $k=0$, and
\begin{equation}\label{ej:F-plane}
	\mathbf F(q)=\varphi(q)^2\in \Pi,\qquad
	q=(x,y,z),
\end{equation}
i.e
\[
	\mathbf F(q)=[(x+y)^2+2(x-z)(y+z)]\, e
	+
	[(y+z)^2-2(x+y)(x-z)]\,e_2+
	[(x-z)^2-2(x+y)(y+z)]\, e_3.
\]
A straightforward calculation yields
\[
	\Delta\,\mathbf F=\mathbf 0,\qquad
	F^i_{xx}+F^i_{yy} +F^i_{zz}=0,\quad
	i=1,2,3.
\]
While, 
\[
	\mathrm{div}\, \mathbf F=-4x+4y+8z,
	\quad
	\mathrm{curl}\,\mathbf F=-4(x+2y+z)e-4(2x+y-z)e_3.
\]
\eejem

\section{Lamellar vector fields}\label{s4}

We consider a more precise description of the vector fields proposed in Theorem \ref{teo:0}. Let us consider a vector field in $\mathbb R^3$, parallel to the plane $\Pi$,
\begin{equation}\label{eq:ej1}
	\mathbf F(x,y,z)={\rm u }(\zeta,\xi,\eta)\,v_2 + {\rm v}(\zeta,\xi,\eta) \,v_3
\end{equation}
Here, we use the linear change of coordinates
\[
	(x,y,z)=   \zeta v_1    + \xi v_2+\eta v_3,
\]
i.e.
\begin{equation}\label{eq:change}
	\left(\begin{array}{c}
	x\\
	y\\
	z
	\end{array}\right)
	=\frac{1}{3}
	\left(\begin{array}{ccc}
		1 & 2 & 0	       \\
		1& -1& \sqrt{3}\\
		1 & -1 &-\sqrt{3}
	\end{array}\right)
	\left(\begin{array}{c}
	\zeta\\
	\xi\\
	\eta
	\end{array}\right).
\end{equation}
From the linear change of coordinates \eqref{eq:change} we obtain,
\begin{equation}\label{eq:change}
	\left(\begin{array}{c}
	F^1\\
	F^2\\
	F^3
	\end{array}\right)
	=\frac{1}{3}
	\left(\begin{array}{ccc}
		1 & 2 & 0	       \\
		1 & -1& \sqrt{3}\\
		1 & -1 &-\sqrt{3}
	\end{array}\right)
	\left(\begin{array}{c}
	0\\
	\mathrm u\\
	\mathrm v
	\end{array}\right).
\end{equation}
The values of the partial directional derivatives $\partial_\alpha F^j$, $i,j=1,2,3,$ $\alpha\in\{\zeta,\xi,\eta\}$ are
\[\begin{array}{rclrclrcl}
	F^1_\zeta&=&{2}\mathrm u_\zeta/3,&
	F^1_\xi&=&{2}\mathrm u_\xi/3,&
	F^1_\eta&=&{2}\mathrm u_\eta/3,
	\\
	F^2_\zeta&=&\left(-\mathrm u_\zeta+\sqrt{3}\mathrm v_\zeta\right)/3,&
	F^2_\xi&=&\left(-\mathrm u_\xi+\sqrt{3}\mathrm v_\xi\right)/3,&
	F^2_\eta&=&\left(-\mathrm u_\eta+\sqrt{3}\mathrm v_\eta\right)/3,
	\\
	F^3_\zeta&=&\left(-\mathrm u_\zeta-\sqrt{3}\mathrm v_\zeta\right)/{3},&
	F^3_\xi&=&\left(-\mathrm u_\xi-\sqrt{3}\mathrm v_\xi\right)/{3},&
	F^3_\eta&=&\left(-\mathrm u_\eta-\sqrt{3}\mathrm v_\eta\right)/{3}.
\end{array}\]
On the other hand, such partial derivatives correspond to the following directional derivatives
\[\begin{aligned}
	F^i_\zeta&:=D_{v_1}F^i(x(\zeta,\xi,\eta),y(\zeta,\xi,\eta),z(\zeta,\xi,\eta)),
	\\
	F^i_\xi&:=D_{v_2}F^i(x(\zeta,\xi,\eta),y(\zeta,\xi,\eta),z(\zeta,\xi,\eta)),
	\\
	F^i_\eta&:=D_{v_3}F^i(x(\zeta,\xi,\eta),y(\zeta,\xi,\eta),z(\zeta,\xi,\eta)).
\end{aligned}\]
Thus, the value of the nine linear variables $F^i_a$, $i\in\{1,2,3\},$ $a\in\{x,y,z\}$, can be obtained by solving the non-degenerate nine equations linear system,
\[
\left(\begin{array}{c}
	F^i_\zeta\\
	F^i_\xi\\
	F^i_\eta
	\end{array}\right)
	=\frac{1}{3}
	\left(\begin{array}{ccc}
		1 & 1 & 1       \\
		2 & -{1}& -{1}\\
		0 & {\sqrt{3}} &-{\sqrt{3}}
	\end{array}\right)
	\left(\begin{array}{c}
	F^i_x\\
	F^i_y\\
	F^i_z
	\end{array}\right).
\]
Whence,
\begin{equation}\label{eq:change2}
	\left(\begin{array}{c}
	F^i_x\\
	F^i_y\\
	F^i_z
	\end{array}\right)
	=
	\left(\begin{array}{ccc}
		1 & 1 & 0	       \\
		1 & -\frac{1}{2}& \frac{\sqrt{3}}{2}\\
		1 & -\frac{1}{2} &-\frac{\sqrt{3}}{2}
	\end{array}\right)\left(\begin{array}{c}
	F^i_\zeta\\
	F^i_\xi\\
	F^i_\eta
	\end{array}\right).
\end{equation}
Hence,
\begin{equation}\label{eq:Fia}
	\begin{aligned}
	F^1_x&={2}(\mathrm u_\zeta+\mathrm u_\xi)/3,
	\\
	F^1_y&={2}\mathrm u_\zeta/3-\mathrm u_\xi/3+\mathrm u_\eta/\sqrt{3},
	\\
	F^1_z&={2}\mathrm u_\zeta/3
	-\mathrm u_\xi/3-\mathrm u_\eta/\sqrt{3},
	\\
	F^2_x&=\left(-\mathrm u_\zeta+\sqrt{3}\mathrm v_\zeta
	-\mathrm u_\xi+\sqrt{3}\mathrm v_\xi\right)/3,
	\\
	F^2_y&=\left(-2\mathrm u_\zeta+2\sqrt{3}\mathrm v_\zeta+
	\mathrm u_\xi-\sqrt{3}\mathrm v_\xi
	-\sqrt{3}\mathrm u_\eta+3\mathrm v_\eta\right)/6,
	\\
	F^2_z&=\left(
	-2\mathrm u_\zeta+2\sqrt{3}\mathrm v_\zeta
	+\mathrm u_\xi-\sqrt{3}\mathrm v_\xi+
	\sqrt{3}\mathrm u_\eta+{3}\mathrm v_\eta\right)/6,
	\\
	F^3_x&=\left(-\mathrm u_\zeta-\sqrt{3}\mathrm v_\zeta-\mathrm u_\xi-\sqrt{3}\mathrm v_\xi\right)/{3},
	\\
	F^3_y&=\left(-2\mathrm u_\zeta-2\sqrt{3}\mathrm v_\zeta+\mathrm u_\xi+\sqrt{3}\mathrm v_\xi
	-\sqrt{3}\mathrm u_\eta-3\mathrm v_\eta\right)/6,
	\\
	F^3_z&=\left(-2\mathrm u_\zeta-2\sqrt{3}\mathrm v_\zeta+
	\mathrm u_\xi+\sqrt{3}\mathrm v_\xi+
	\sqrt{3}\mathrm u_\eta-{3}\mathrm v_\eta\right)/6.	
\end{aligned}\end{equation}

$
	\mathbf F
$ is $\varphi\mathbb A$-differentiable if by substitution of \eqref{eq:Fia} in CR relations \eqref{eq:system2} we obtain four linearly independent relations in six variables, $\mathrm u_\zeta,\mathrm u_\xi,\mathrm u_\eta,\mathrm v_\zeta,\mathrm v_\xi,\mathrm v_\eta$. Namely,
\begin{equation}\label{eq:CR-uv}
\begin{aligned}
	2\mathrm u_\zeta-\frac{1}{2}\mathrm u_\xi-\frac{\sqrt{3}}{2}\mathrm u_\eta
	+\frac{\sqrt{3}}{2}\mathrm v_\xi	-\frac{1}{2}\mathrm v_\eta
	&=0,
	\\
	-\mathrm u_\zeta+\frac{1}{2}\mathrm u_\xi+\frac{\sqrt{3}}{2}\mathrm u_\eta
	-{\sqrt{3}}\mathrm v_\zeta
	-\frac{\sqrt{3}}{2}\mathrm v_\xi-\frac{1}{2}\mathrm v_\eta
	&=0,
	\\
	-\mathrm u_\zeta-\mathrm u_\xi+\frac{1}{2}\mathrm v_\zeta
	-\frac{\sqrt{3}}{2}\mathrm v_\xi+
	\mathrm v_\eta
	&=0,
	\\
	\frac{3}{2}\mathrm u_\xi-\frac{\sqrt{3}}{2}\mathrm u_\eta
	-\frac{2}{3}\mathrm v_\zeta
	-\frac{1}{2\sqrt{3}}\mathrm v_\xi-\frac{1}{2}\mathrm v_\eta
	&=0.
\end{aligned}\end{equation}

\bcor
A vector field \eqref{eq:ej1} satisfying CR equations \eqref{eq:CR-uv}. Then $\Delta\mathrm u=0$ in the open set $\mathcal U$, i.e. the component $\mathrm u$ satisfies
\begin{equation}
\begin{aligned}
\mathrm u_{xx}+\mathrm u_{yy}+\mathrm u_{zz}&=0,
\\
\mathrm u_{xxx}+\mathrm u_{yyx}+\mathrm u_{zzx}&=0,
\\
\mathrm u_{xxy}+\mathrm u_{yyy}+\mathrm u_{zzy}&=0,
\\
\mathrm u_{xxz}+\mathrm u_{yyz}+\mathrm u_{zzz}&=0.
\end{aligned}
\end{equation}
The component $\mathrm v$ of $\mathbf F$ also satisfies such equations.
\ecor

\bthm\label{teo:1} With the same hypothesis as in Theorem \ref{teo:0}, the vector field $\mathbf V:\mathcal U\rightarrow\mathbb A$ defined using \eqref{eq:V(F)} as
\[
	\mathbf V =V(\mathbf F)= ( F^3-F^2)e_1+
			( F^3- F^1)e_2+
			( F^2- F^1)e_3\in V(\Pi)
\]
satisfies
\[
	\mathrm{div}\, \mathbf V=0,\quad
	\mathrm{curl}\, \mathbf V={\bf 0}.
\]
Moreover, $\mathbf V$ is also harmonic,
\[
	\Delta \mathbf V=\mathbf 0.
\]
\ethm

\noindent\textbf{Proof of Theorem \ref{teo:1}.} 
The following equations can be deduced fro CR relations,
\begin{equation}\label{eq:VF}
	 F^2_x - F^1_x
	 = 
	 F^3_z -F_z^2,
	\qquad
	F^2_y - F^1_y
	 = 
	 F^3_z -  F_z^1,
	\qquad
	 F^3_x -  F^1_x
	 = 
	F^3_y -  F_y^2.
\end{equation}
For $\mathbf V= ( F^3-F^2)e_1+
			( F^3- F^1)e_2+
			( F^2- F^1)e_3$, \eqref{eq:VF} imply,
\[
	 V^3_x= V^1_z,
	 \quad
	V^3_y= V^2_z,
	\quad
	V^2_x=V^1_y
	 .	
\]
or $\mathrm{curl}\,\mathbf V=0$. Similarly,
\[
	2\left[ V^1_x+V^2_y+ V^3_z\right]
	=
	 2\left[( F^3_x- F^2_x)+
			( F^3_y- F^1_y)+
			( F^2_z- F^1_z)\right].
\]

And from \eqref{eq:system2} we get
\[\begin{aligned}
	2\,\mathrm{div}\, \mathbf V
	&=
	( F^1_y-F^2_y+ F^1_z- F^3_z)
	+
	( F^2_x- F^1_x+ F^3_z- F^2_z)
	+(F^2_y-F^3_y+ F^1_x- F^3_x)
	\\
	&=
	F^1_y+ F^1_z+ F^2_x- F^2_z- F^3_y-F^3_x
	\\
	&=-( F^3_y-F^1_y)-( F^3_x- F^2_x)-( F^2_z- F^1_z)
	\\
	&=-{\rm div}\, \mathbf V
\end{aligned}
\]
which implies that $\mathrm{div}\, \mathbf V=0$.
\hfill $\Box$

\brem\label{rmk:22}
$\mathbf V$ is no longer $\varphi\mathbb A$-differentiable. In fact, vector $\mathbf V$ can be constructed using solely vectors $\mathbf F$ parallel to $\Pi$. Since $V(\mathbb R^3)= V(\Pi)$ then the vector field $V(\mathbf F)$ is a flat vector field parallel to the plane $V(\Pi)$. In its turn, if ${\rm curl}\,\mathbf V=\mathbf 0$, then in a simply connected domain in the plane $V(\Pi)$, the flat vector field $V(\mathbf  F)$ would be a gradient-like vector field in the plane $V(\Pi)$. Since $\mathbf V$ is also divergence-free, then the potential function in its turn would be harmonic along $V(\Pi)$. 
\erem

\bejem
When we consider a polynomial vector field $\mathbf F(q)=c_0+c_1\centerdot \varphi(q)+\dots+c_2\centerdot \varphi(q)^2$ with $c_i\in\mathbb A$, then we have $\mathbf V=V(\mathbf F)$ irrotational and incompressible. Take for instance $\mathbf F$ as in  \eqref{ej:F-plane}, then
\[\begin{aligned}
\mathbf V&=
\left\{[(x-z)^2-2(x+y)(y+z)]-[(y+z)^2-2(x+y)(x-z)]\right\}\, e\\
 &+ 
 \left\{
	[(x-z)^2-2(x+y)(y+z)]-[(y+z)^2-2(x+y)(x-z)]\right\}\,e_2\\
&+\left\{
	[(y+z)^2-2(x+y)(x-z)]-[(x+y)^2+2(x-z)(y-z)]\right\}\, e_3
	.
\end{aligned}\]
satisfies the conclusions of Theorem \ref{teo:1}. Indeed, a straightforward calculation yields
\[
	\Delta\,\mathbf V
	=\mathbf 0,\qquad
	V^i_{xx}+V^i_{yy} +V^i_{zz}=0,\quad
	i=1,2,3,
\]
and, 
\[
	\mathrm{div}\, \mathbf V=0,
	\quad
	\mathrm{curl}\,\mathbf V=\mathbf 0.
\]
which is equivalent to the system of equations \eqref{eq:desglose}.
\eejem

The $\varphi\mathbb A$-differentiable vector fields $\mathbf F$ described in \eqref{eq:ej1} have linear first integral $H(x,y,z)=x+y+z$. The corresponding harmonic vector fields $\mathbf V=V(\mathbf F)$ have $ H_1(x,y,z)=x-y+z$ as linear first integral.

When we consider a $\Pi$-parallel vector field $\mathbf F$ as in \eqref{eq:ej1} we obtain a $V(\Pi)$-parallel vector field,
\[
	\mathbf V=
	\mathrm u \,w_2+\mathrm v\, w_3
=	
	\mathrm u \,V(v_2)+\mathrm v\, V(v_3)
	=
	-\frac{2}{\sqrt{3}}\mathrm v\, e
	+\left(\mathrm u-\frac{1}{\sqrt{3}}\mathrm v\right)\,e_2
	+\left(\mathrm u+\frac{1}{\sqrt{3}}\mathrm v\right)\,e_3,
\]
where $(\mathrm u,\mathrm v)$ satisfy the induced CR relations \eqref{eq:CR-uv}. Then, relations \eqref{eq:solenoid-vf} become the linear system \eqref{eq:desglose} below. \begin{equation}\label{eq:desglose}
\begin{aligned}
	\mathrm u_y+\mathrm u_z-\frac{2}{\sqrt{3}}\mathrm v_x-\frac{1}{\sqrt{3}}\mathrm v_y
	+\frac{1}{\sqrt{3}}\mathrm v_z
	&=0,
	\\
	\mathrm u_y-\mathrm u_z+\frac{1}{\sqrt{3}}\mathrm v_y
	+\frac{1}{\sqrt{3}}\mathrm v_z
	&=0,
	\\
	\mathrm u_x+\frac{1}{\sqrt{3}}\mathrm v_x+\frac{2}{\sqrt{3}}\mathrm v_z
	&=0,
	\\
	\mathrm u_x-\frac{1}{\sqrt{3}}\mathrm v_x+\frac{2}{\sqrt{3}}\mathrm v_y
	&=0.
\end{aligned}\end{equation}
System \eqref{eq:desglose} consists of 4 independent equations in 6 variables given by partial derivatives of 
\[
	\mathrm{u}(\zeta(x,y,z),\xi(x,y,z),\eta(x,y,z)),
	\qquad
	\mathrm{v}(\zeta(x,y,z),\xi(x,y,z),\eta(x,y,z)),
\]
with respect to $x,y,z$, respectively. Thus, there are 2 directional derivatives for $\mathrm u$ and/or $\mathrm v$ that can be chosen freely, while the remaining 4 partial derivatives are constrained by \eqref{eq:desglose}. In particular, the following couple of CR equations on directional derivatives along the basis, $w_2=V(v_2),$ and $w_3=V(v_3)$ defined in Lemma \ref{lma:3}, are implied by \eqref{eq:desglose}
\begin{equation}\label{eq:2D-CR}
	D_{w_2}\mathrm u=-D_{w_3}\mathrm v,\qquad
	D_{w_3}\mathrm u=\frac{1}{3}\,D_{w_2}\mathrm v,
\end{equation}
which accurately describe lamellar vector fields $\mathbf V$ as $1$-parameter couples of functions $\mathrm u^\zeta=\mathrm u(\cdot,\cdot,\zeta),$ and $\mathrm v^\zeta=\mathrm v(\cdot,\cdot,\zeta)$ depending differentiably on $\zeta\in\mathbb R$ and solving \eqref{eq:2D-CR}.


\end{document}